\documentclass[11pt]{article}

\usepackage{siunitx}
\usepackage{booktabs}
\usepackage{amsmath} 
\usepackage{color}
\usepackage{mathtools}
\usepackage{amsfonts}
\usepackage{amsthm}
\usepackage{float}
\usepackage{color}

\DeclarePairedDelimiter\abs{\lvert}{\rvert}

\title{The Runge Example for Interpolation and\\Wilkinson's Examples for Rootfinding}
\author{Robert M. Corless and Leili Rafiee Sevyeri}
\date{}

\usepackage{xfrac}

\usepackage{csquotes}

\usepackage{listings}
\theoremstyle{definition}

\newtheorem*{remark}{Remark}
\renewcommand{\emph}[1]{\textsl{#1}}

\newtheorem{theorem}{Theorem}

\newcommand{\ma}{\textsc{MATLAB}}

\begin{document}

\clearpage\maketitle
\thispagestyle{empty}

\noindent\rule{\textwidth}{0.5mm}
\begin{abstract}
We look at two classical examples in the theory of numerical analysis, namely the Runge example for interpolation and Wilkinson's example (actually two examples) for rootfinding. We use the modern theory of backward error analysis and conditioning, as instigated and popularized by Wilkinson, but refined by Farouki and Rajan. By this means, we arrive at a satisfactory explanation of the puzzling phenomena encountered by students when they try to fit polynomials to numerical data, or when they try to use numerical rootfinding to find polynomial zeros. Computer algebra, with its controlled, arbitrary precision, plays an important didactic role.
\end{abstract}

The function $y = \frac{1}{1 + 25x^2}$ on $-1 \leq x \leq 1$, called the Runge example\footnote{Carl David Tolm\'e Runge ($30$ August 
$1856-3$ January $1927$) was a German mathematician, physicist, and spectroscopist.}, is used in many numerical analysis textbooks to show why high-degree polynomial interpolation \textsl{using equally-spaced nodes} is bad.
Unfortunately, most textbooks omit or downplay the crucial qualification ``using equally-spaced nodes'' and thereby leave the \textsl{false} impression that high-degree interpolation is \textsl{always} bad.

Similarly, Wilkinson's first example polynomial\footnote{James Hardy Wilkinson. Born: $27$ September $1919$ in Strood, Kent, England- Died: $5$ October $1986$ in Teddington, Middlesex, England. He worked with Turing in $1946$. He won the Chauvenet Prize in $1970$ for mathematical exposition for his paper ``The Perfidious Polynomial"~\cite{wilkinson}. His book, ``The Algebraic Eigenvalue Problem" was foundational for the field of numerical linear algebra.} 
\begin{equation}
    p_{20}(x) = \prod_{k=1}^{20}(x-k)
\end{equation}
is widely discussed as an example ---maybe \textsl{the} canonical example--- of polynomial perfidy, this time for rootfinding, not interpolation.
As we will discuss below, both examples are better explained using the theory of \textsl{conditioning}, as developed for instance by Farouki and Rajan~\cite{farouki1987numerical}. One considers a polynomial
\begin{equation}\label{eq2}
p(x) = \sum_{k=0}^n c_k\phi_k(x)
\end{equation}
expressed in some polynomial basis $\{\phi_k(x)\}_{k=0}^n$ for polynomials of degree at most $n$.
The usual monomial basis $\phi_k(x) = x^k$ is the most common, but by no means best for all purposes. Farouki and Goodman~\cite{farouki1996optimal} point out that the Bernstein basis $\phi_k(x) = {n \choose k }x^k(1-x)^{n-k}$ has the best conditioning in general, out of all polynomial bases satisfying $\phi_k(x)\geq 0$ on the interval $0\leq x\leq 1$. Surprisingly, Corless and Watt~\cite{corless2004bernstein} show Lagrange bases can be better; see also J.M. Carnicer and Y. Khiar~\cite{Carnicer}.\\

We now discuss Farouki and Rajan's formulation. The idea is that one investigates the effects of small relative changes to the coefficients $c_k$, such as might arise from data error or perhaps approximation or computational error.
The model is
\begin{equation}\label{eq3}
    p(x) + \Delta p(x) = \sum_{k=0}^n c_k (1 + \delta_k) \phi_k(x)
\end{equation}
where each $| \delta_k| \le \varepsilon$, usually taken to be small.
For instance, if $\varepsilon = 0.005$, each coefficient can be in error by no more than $0.5\%$.
In particular, zero coefficients are not allowed to be disturbed at all.

Then,
\begin{align}
    | \Delta p(x) | &= \left | \sum_{k=0}^n c_k (1 + \delta_k) \phi_k(x) - \sum_{k=0}^n c_k \phi_k(x) \right |\\
    &= \left | \sum_{k=0}^n c_k \delta_k \phi_k(x) \right |. \label{5}
\end{align}
By the triangle inequality, this is

\begin{align}
    &\le \sum_{k=0}^n \left | c_k \delta_k \phi_k(x) \right | \\
    &\le \left ( \sum_{k=0}^n |c_k| |\phi_k(x)| \right ) \max_{0 \le k \le n} | \delta_k|\\
    &\le B(x) \cdot \varepsilon
\end{align}
where
\begin{equation}\label{9}
    B(x) = \sum_{k=0}^n |c_k| |\phi_k(x)|
\end{equation}
serves, for each $x$, as a \textsl{condition number} for polynomial evaluation.\\
This is to be contrasted with the definition of ``evaluation condition number" that arises when thinking of error $\Delta x$ in the input: If $x$ changes to $x+\Delta x$, then $y=f(x)$ changes to $y+\Delta y$ where calculus tells us that 
\begin{equation}
\dfrac{\Delta y}{y}\doteq \dfrac{xf'(x)}{f(x)}\cdot\dfrac{\Delta x}{x}\>.
\end{equation}
Here $C=xf'(x)/f(x)$ is the condition number, and instead of $B(x)$; but $B$ and $C$ measure the responses to different types of error (coefficients and input). We look at $B$, here.\\
\begin{remark}
There are many theorems\footnote{See for instance those cited in chapter $2$ of~\cite{corless2013graduate}.} in numerical analysis that say, effectively, that when evaluating equation \eqref{eq2} in IEEE floating-point arithmetic, that the computed result is exactly of the form \eqref{eq3} for $|\delta _k|<K.\mu$ where $K$ is a modest constant and $\mu$ is the unit roundoff --- in double precision, $2^{-53}\doteq 10^{-16}$. This is one motivation to study the effects of such perturbations, but there are others.
\end{remark}


\subsubsection*{The Runge Example}
\begin{figure}[!h]
    \centering
    \includegraphics[width=0.5\textwidth]{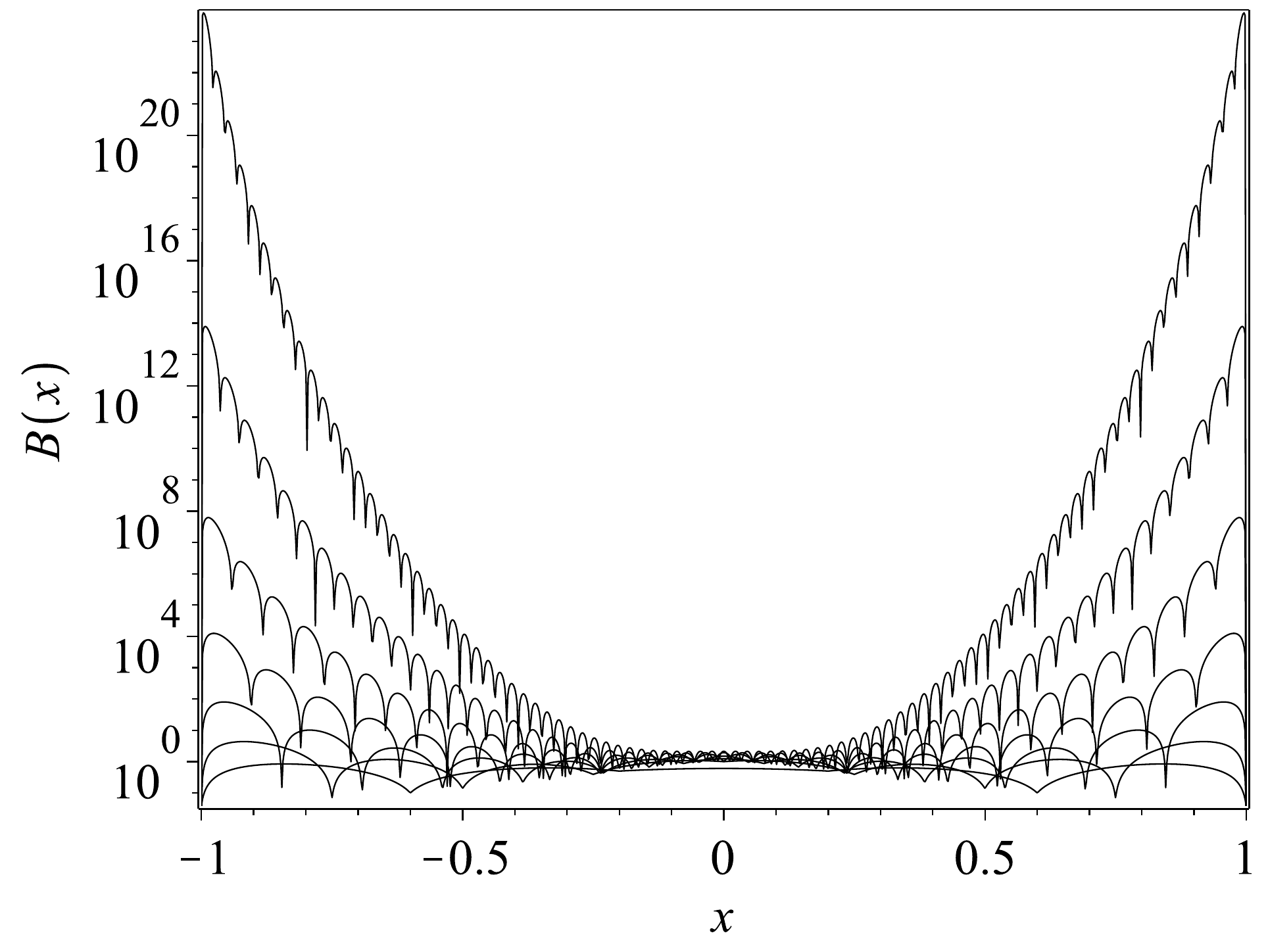}
    \caption{The condition number of the Runge example on equally-spaced nodes with degrees $n=5, 8, 13, 21, 34, 55, \text{and } 89$.}
    \label{Figure1}
\end{figure}

If the equally-spaced nodes $x_k = -1 + 2k/n$, $k=0\ldots n$ are used to interpolate a function with a single polynomial of degree at most $n$, and the basis functions
\begin{equation}
    \ell_k(x) = \frac{\prod\limits_{\substack{{j=0}\\{j \ne k}}}^n (x - x_j)}{\prod\limits_{\substack{{j=0}\\{j \ne k}}}^n (x_k - x_j)}
\end{equation}
which are the Lagrange interpolation basis functions, are used, then the coefficients are the values for the Runge function:
\begin{equation}
    y_k = f(x_k) = \frac{1}{1 + 25x_k^2} \>.
\end{equation}
Then the condition number of the interpolant is,
\begin{equation}
    B(x) = \sum_{k=0}^n \frac{1}{1 + 25x_k^2} | \ell_k(x) | \>.
\end{equation}
Choosing $n = 5, 8, 13, 21, 34, 55, \text{and } 89$, we plot $B(x)$ on a logarithmic vertical scale for $-1 \le x \le 1$.
The result is in Figure \ref{Figure1}.
The Maple code used to generate that figure is as follows (similar code using \ma\ can be provided, but the Maple code below avoids numerical issues in the construction of $B(x)$ by using exact rational arithmetic).
We see that the maximum values for $B(x)$ occur near $x = \pm 1$, and that for any $n$ there is an interval over which $B(x)$ is small.

\begin{lstlisting}
Digits := 15:
Ns := [seq(combinat[fibonacci](k), k=5..11)]:
f := x -> 1/(1 + 25*x^2):

for N in Ns do
    tau := [seq(-1 + 2*k/N, k=0..N)]:
    rho := [seq(y[k], k=0..N)]:

    p := CurveFitting[PolynomialInterpolation](tau, rho, z, form=Lagrange):

    B := map(abs, p):
    BRunge := eval(B, [seq(y[k] = f(tau[k+1]), k=0..N)]):

    pl[N] := plots[logplot](BRunge, z=-1..1, color=black):

end do:

plots[display]([seq(pl[N], N in Ns)]);
\end{lstlisting}
This experiment well illustrates that interpolation of the Runge example function on equally-spaced nodes is a bad idea. But, really, it is the nodes that are bad.\\
\begin{displayquote}
``Generations of textbooks have warned readers that polynomial interpolation is dangerous. In fact, if the interpolation points are clustered and a stable algorithm is used, it is bulletproof.''
\\[5pt]
\rightline{{\rm --- L.N. Trefethen,~\cite{Trefethen}}}
\end{displayquote}
For a full explanation of the Runge phenomenon, see chapter $3$ of~\cite{Trefethen}.

\subsubsection*{The Runge Example with Chebyshev Nodes}
\begin{figure}[H]
    \centering
    \includegraphics[width=0.5\textwidth]{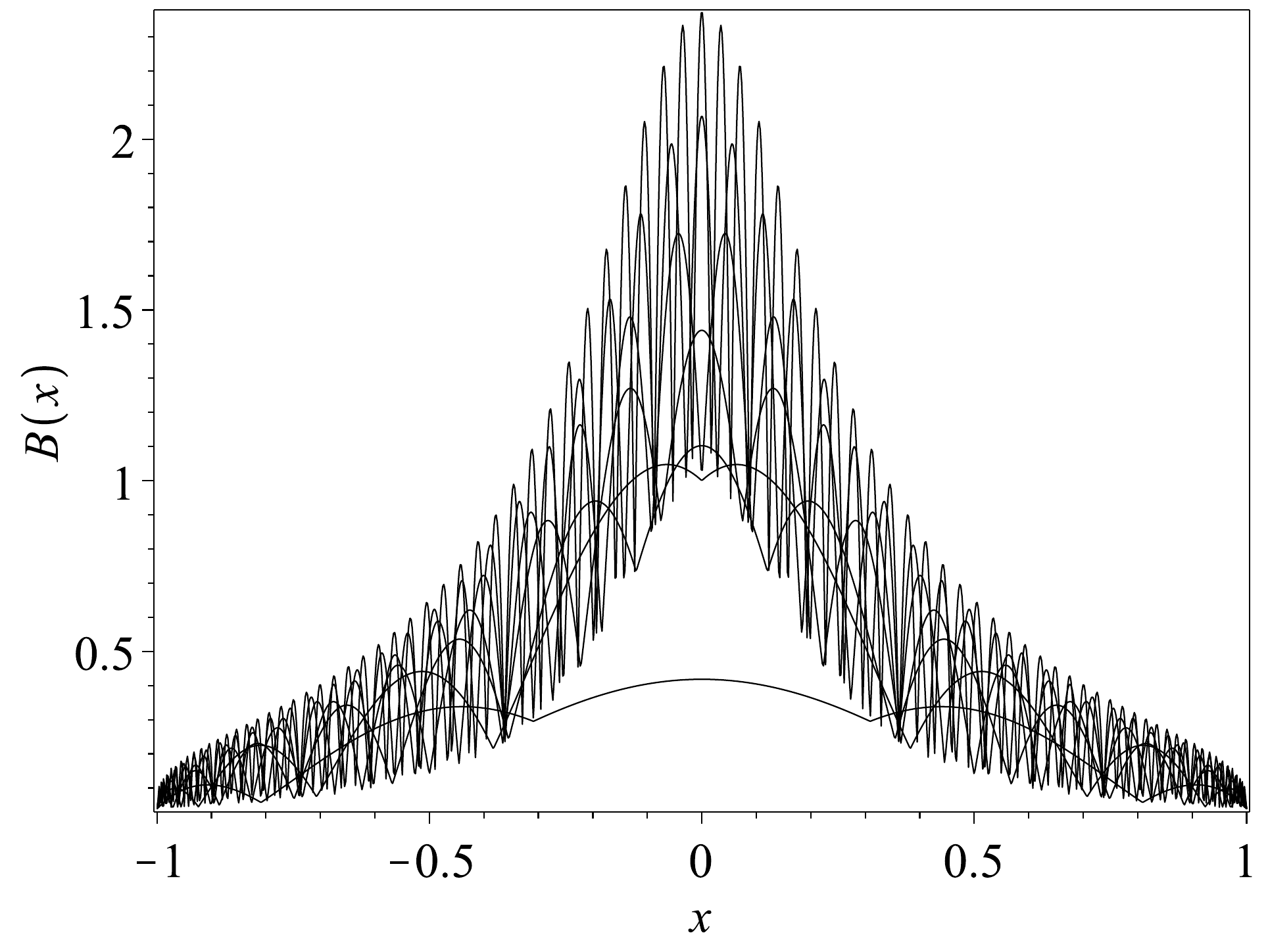}
    \caption{The Runge example with Chebyshev nodes with degrees $n=5, 8, 13, 21, 34, 55, \text{and } 89$.}
    \label{plot2}
\end{figure}

If instead we use $x_k = \cos(\pi k/n)$, replacing the line
\begin{lstlisting}[frame=none]
    tau := [seq(-1 + 2*k/N, k=0..N)];
\end{lstlisting}
with
\begin{lstlisting}[frame=none]
 tau := [seq(evalf[2*Digits](cos(Pi*k/N)), k=0..N)];
\end{lstlisting}
then $B(x)$ climbs no higher than about 2.
Indeed we can replace \texttt{plots[logplot]} by just \texttt{plot}. See Figure \ref{plot2}.\\
This is an improvement, for $n=89$, of about a factor of $10^{22}$.
For a detailed exposition of why this works, and when, see~\cite{corless2013graduate} and the Chebfun project at www.chebfun.org.

\subsubsection*{Concluding remarks on the Runge example}
An instructor of numerical analysis has to walk a tightrope: the students need to be taught caution (maybe bordering on paranoia) but they also need to learn when to trust their results. Learning to assess the sensitivity of their expression (as programmed) to realistic changes in data values is an important objective. It is true that people (not just people who are students) don't want correct but complicated answers, preferring   simple answers that they don't have to worry about. The Runge example is a very clear case where these ideas can be usefully and thoroughly explored. \\
On can go further and replace $B(x)$ by its upper bound in terms of the Lebesgue function $B(x)\leq L(x)||c||_{\infty}$ where $L(x)=\sum_{k=0}^{n}|\phi_k(x)|^k$. This more general analysis is useful, as in~\cite{Trefethen} but loses in our opinion the chance to make a special retrospective diagnostic of the problem at hand. Moreover, there are cases where $B(x)\ll L(x)||c||_{\infty}$ and this overestimation could lead to the wrong conclusion.\\
The bad behavior of the Runge example shows up in other ways, notably in the ill-conditioning of the Vandermonde matrix on those nodes. But the Vandermonde matrix is ill-conditioned on Chebyshev nodes, too~\cite{Beckermann}; so that can't be the whole story. The explanation offered here seems more apt.

\subsubsection*{Wilkinson's First Polynomial}
Let us now consider rootfinding. Suppose $r$ is a simple zero of $p(x)$: That is, $\,p'(x) \ne 0$ and
\begin{equation}
    0 = p(r) = \sum_{k=0}^n c_k \phi_k(r) \>.
\end{equation}
Suppose $r + \Delta r$ is the corresponding zero of $p + \Delta p$.
This really only makes sense if $\Delta p$ is sufficiently small. Otherwise, the roots get mixed up.
Then
\begin{align}
    0 &= (p + \Delta p)(r + \Delta r) = p(r + \Delta r) + \Delta p(r + \Delta r)\\
    & \approx p(r) + p'(r)\Delta r + \Delta p(r) + \mathcal{O}(\Delta^2)
\end{align}
to first order; since $p(r) = 0$ also we have 
\begin{equation}
    p'(r)\Delta r \approx -\Delta p(r)
\end{equation}
or 
\begin{equation}
    | \Delta r| \approx \left | \frac{-\Delta p(r)}{p'(r)} \right | \le \frac{B(r) \cdot \varepsilon}{|p'(r)|}
\end{equation}
where
\begin{equation}
    B(r) = \sum_{k=0}^n |c_k| |\phi_k(r)|
\end{equation}
as before is the condition number. For nonzero roots, the number \\
$A(r)=\abs*{\frac{rB(r)}{p'(r)}}$ has $\abs*{\frac{\Delta r}{r}}\leq A(r)\varepsilon$ giving a kind of mixed relative/absolute conditioning.
This analysis can be made more rigorous by using ``pseudozeros'' as follows. Define, for given $w_k\geq 0$ not all zero,
\begin{equation}\label{19}
\Lambda_\varepsilon(p) := \{z\, ; \,\exists\,\Delta c_k  \text{ with } |\Delta c_k |\leq w_k \varepsilon \text{ and } \sum_{k=0}^n (c_k + \Delta c_k) \phi_k(z) = 0\}\>.
\end{equation}
Normally, we take $w_k=|c_k|$ in which case we may write $\Delta c_k=c_k\delta_k$.\\
This is the set of all complex numbers that are zeros of ``nearby" polynomials---nearby in the sense that we allow the coefficients to change. This definition is inconvenient to work with. Luckily, there is a useful theorem, which can be found, for instance, in~\cite[Theorem $5.3$]{corless2013graduate}, also see~\cite{Green} and~\cite{amir}.\\

\begin{figure}[H]
    \centering
    \includegraphics[width=0.5\textwidth]{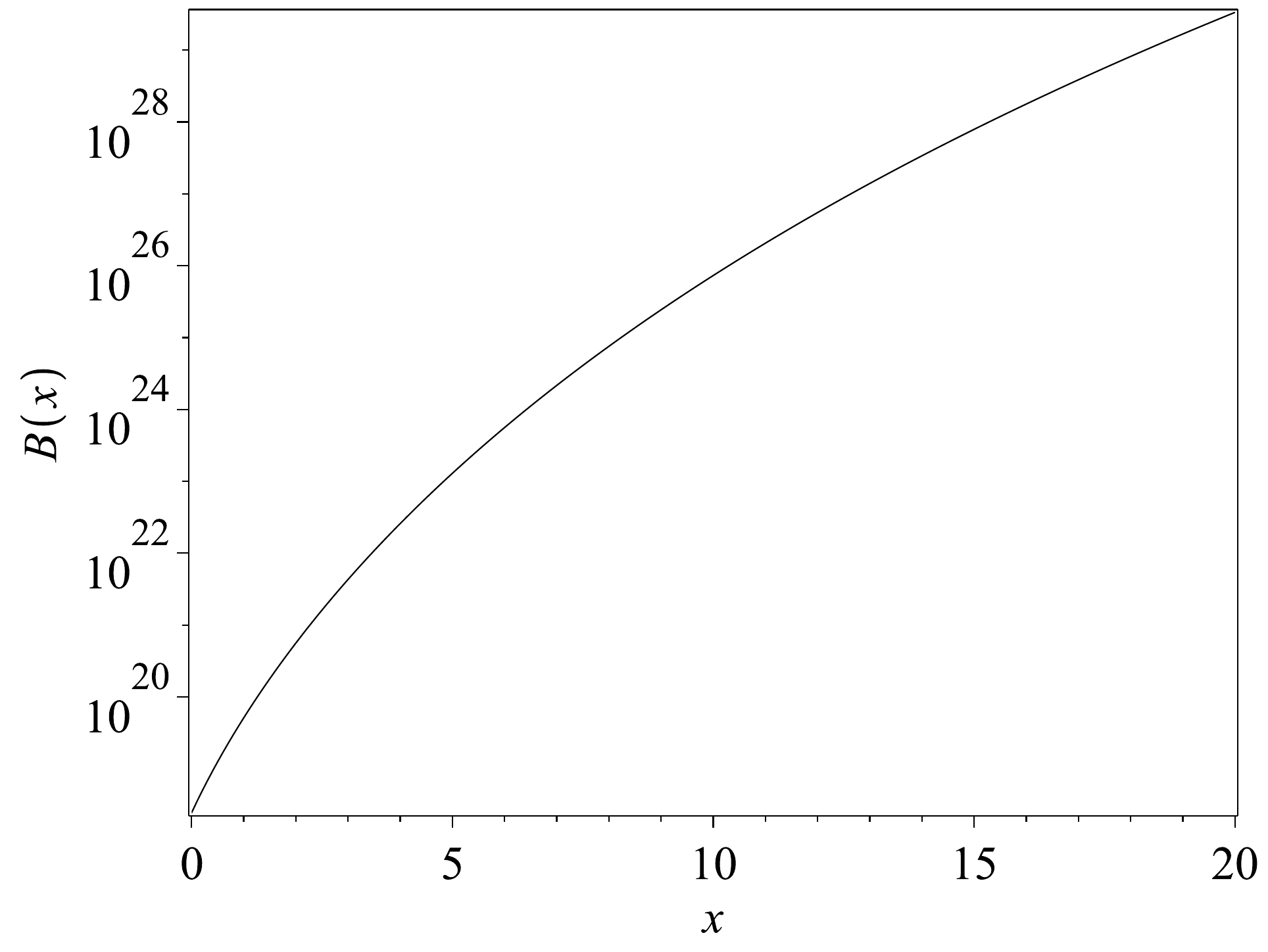}
    \caption{ The condition number of Wilkinson's first polynomial ($N=20$).}
    \label{plotw}
\end{figure}

\begin{figure}[H]
    \centering
    \includegraphics[width=0.4\textwidth]{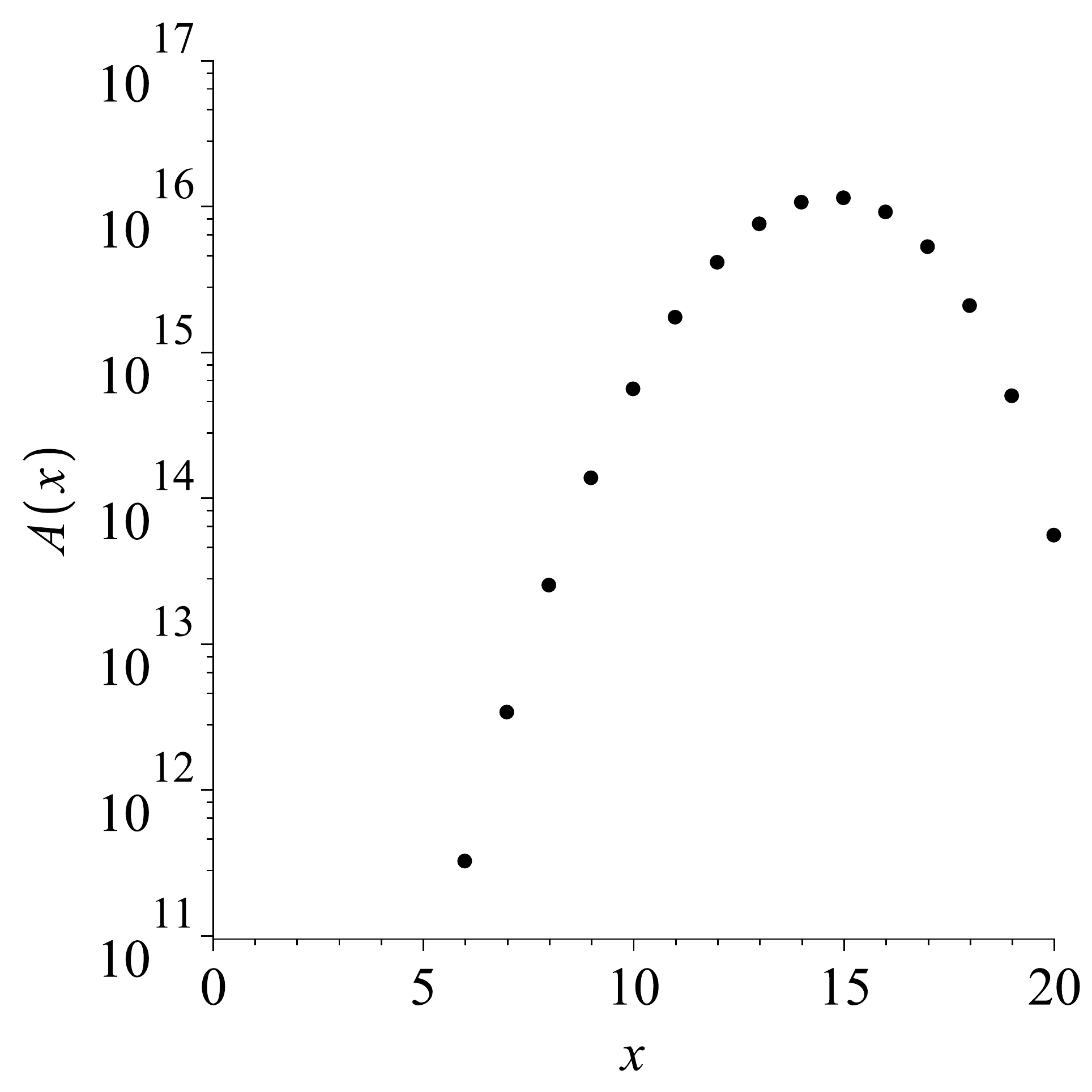}
    \caption{The condition number for rootfinding.}
    \label{ampli}
\end{figure}

\begin{theorem}
Given weights $w_k\geq 0$, not all zero, and a basis $\phi_k(z)$, define the weighted $\varepsilon$-pseudozero set of $p(z)$ as in equation \eqref{19}.
Suppose also that 
$$\delta p(z)=\sum_{k=0}^{n}\Delta c_k \phi_k(z).$$
Moreover, let 
$$B(\lambda)=\sum_{k=0}^{n}w_k |\phi_k(\lambda)|.$$
Then the pseudozero set of $p(z)$ may be alternatively characterized as
\begin{equation}
    \Lambda_\varepsilon(p)=\{
    z\, ; \, |p(z)| \le B(z) \cdot \varepsilon
    \} = \lbrace z\, ; \, \abs*{\frac{zp(z)}{p'(z)}}\leq \abs*{\frac{zB(z)}{p'(z)}}\varepsilon \rbrace 
\end{equation}
\end{theorem}
This is again a condition number; the same one, as in equation \eqref{9} if $w_k = |c_k|$.\\
Wilkinson's first polynomial is, with $N = 20$,
\begin{align}
    W_N(x) &= \prod_{k=1}^N (x - k)\\
    &= (x-1)(x-2)(x-3)\cdots(x-N) \>.
\end{align}
In this form, it is ``bulletproof''.
However, if we are so foolish as to expand it into its expression in the monomial basis, namely,
\begin{equation}
    W_N(x) = x^N - \frac{1}{2}N(N+1)x^{N-1} + \cdots + (-1)^N \cdot N!
\end{equation}
(for $N = 20$ this is $x^{20} - 210x^{19} + \cdots + (20!)$) and in this basis, $\phi_k = x^k$, the condition number for evaluation 
\begin{equation}
    B_N(x) = |x|^N + \frac{1}{2}N(N+1)|x|^{N-1} + \cdots + |N!| \>.
\end{equation}
is very large. See Figure \ref{plotw}.
When we plot the condition number for root finding, $A(r)=\abs*{\frac{r B_N(r)}{W'_N(r)}}$, we find that for $N=20$ (Wilkinson's original choice), the maximum value occurs at $r = 16$ and $r B_{20}(r)/| W'_{20}(r)| \approx 10^{16}$. See Figure \ref{ampli}.
Working in single precision would give no figures of accuracy; double precision ($u \approx 10^{-16}$) also does not guarantee any accuracy.
For $N = 30$ we find $\frac{r B_{30}(r)}{|W'_{30}(r)|} > 10^{21}$ sometimes; for $N = 40$ it's $10^{28}$.
Working with the monomial basis for this polynomial is surprisingly difficult.
Wilkinson himself was surprised; the polynomial was intended to be a simple test problem for his program for the ACE computer.
His investigations led to the modern theory of conditioning~\cite{wilkinson}.

However, there's something a little unfair about the scaling: the interval $0 \le x \le 20$ when taken to the twentieth power covers quite a range of values.
One wonders if matters can be improved by a simple change of variable.

\subsubsection*{The Scaled Wilkinson Polynomial}
\begin{figure}[H]
    \centering
    \includegraphics[width=0.4\textwidth]{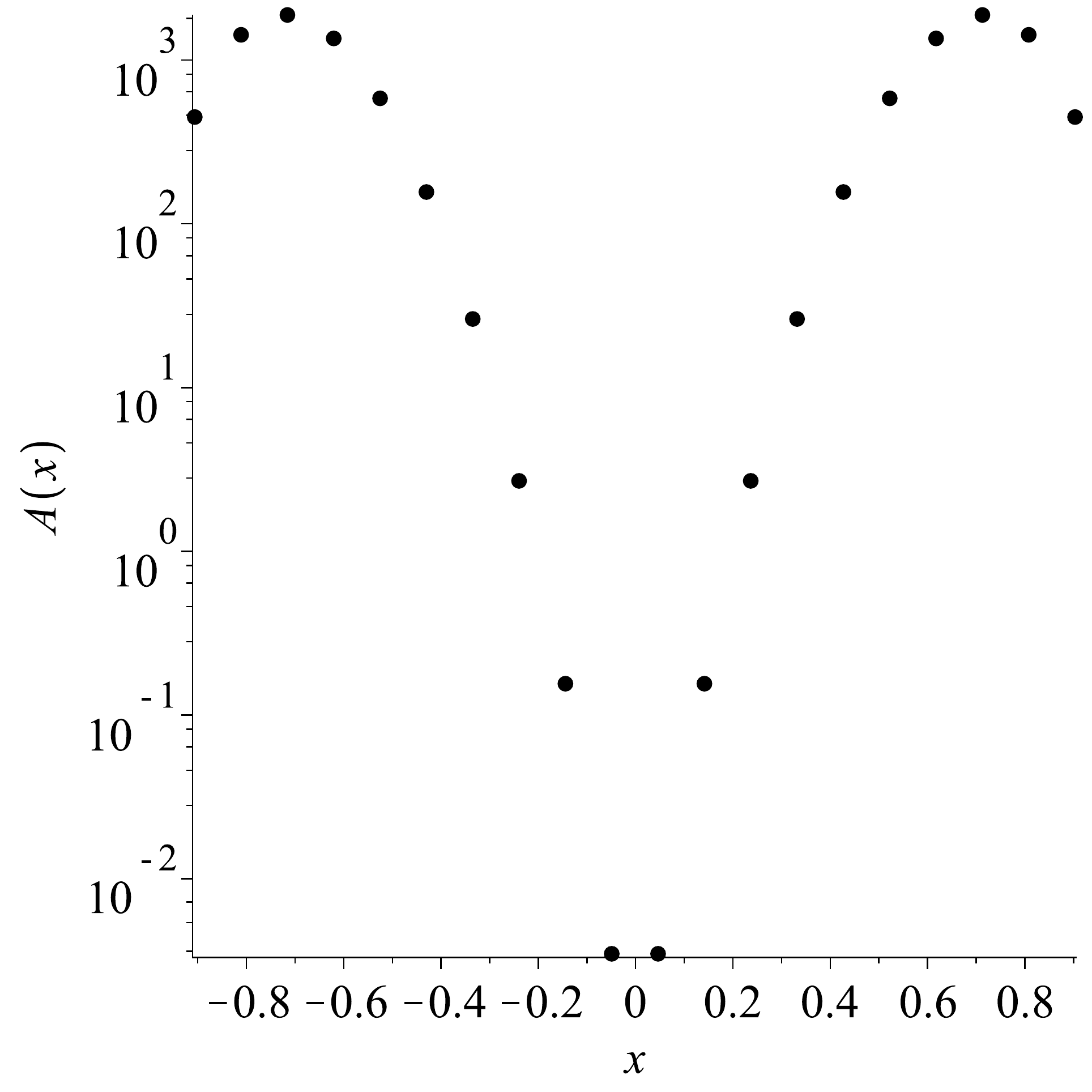}
    \caption{The condition number for the scaled Wilkinson polynomial, $A(r)=\abs*{\frac{r B_N(r)}{W'_N(r)}}$.}
    \label{scaledwil1}
\end{figure}

\begin{figure}[H]
    \centering
    \includegraphics[scale=0.3]{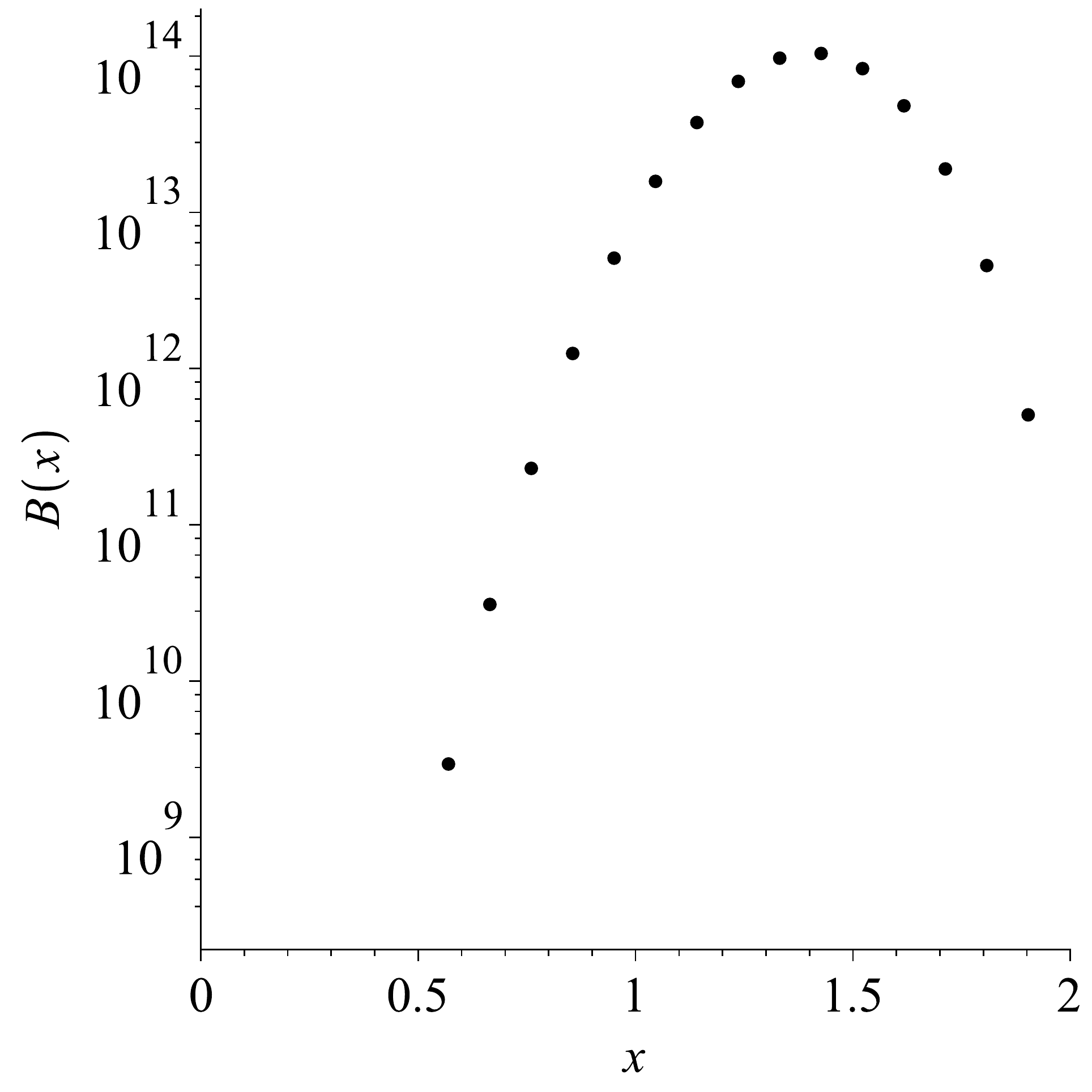}
    \caption{The condition number of scaled Wilkinson polynomial.}
    \label{scaledwil2}
\end{figure}

\begin{figure}[H]
    \centering
    \includegraphics[width=0.5\textwidth]{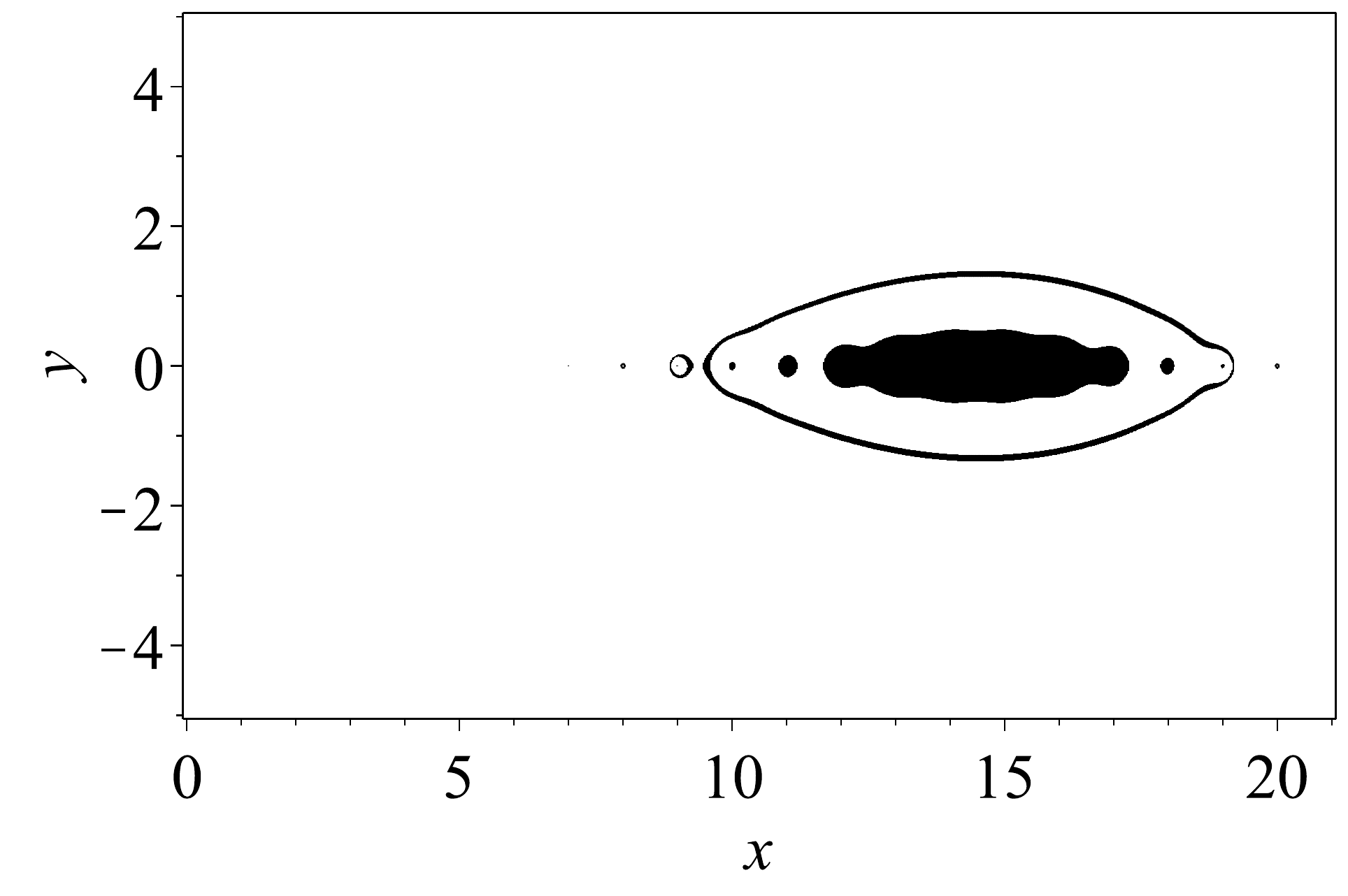}
    \caption{The pseudozeros of $W_N(x)$. The contour levels are $10^{-14}$ and $10^{-18}$. The interior is blacked out because contours are difficult to draw at such sizes, in floating point arithmetic.}
    \label{psedouzeros}
\end{figure}

If we move the roots $1, 2, 3, \ldots, 20$ to the roots $-1 + 2k/21$, $k = 1\ldots30$, then they become symmetrically placed in $-1 < x < 1$, and this improves matters quite dramatically, as we can see in Figure \ref{scaledwil1}.\\
The condition number $10^{13}$ becomes just $10^3$, and we have to go to $N = 60$ to get condition numbers as high as $10^{13}$.
The scaling seems to matter. However, nearly all of the improvement comes from the symmetry; $W_N$ will be even if $N$ is even, and odd if $N$ is odd, and this means half the coefficients are zero and therefore not subject to (relative) perturbation.\\
If instead we scale to the interval $[0,2]$ we have a different story: for roots $2 - 2k/21$ the condition number $B(x)$ reaches nearly the same heights as it did on $0 \le x \le 20$. See Figure~\ref{scaledwil2}. Similarly if we use $[0,1]$. Thus we conclude that symmetry matters.\\
See Figure~\ref{psedouzeros} for the pseudozeros of $W_N(x)$, where the contour levels are $10^{-14}$ and $10^{-18}$. The roots are visibly changed by extremely tiny perturbations. 

\subsubsection*{Wilkinson's Second Example Polynomial}
\begin{figure}[H]
    \centering
    \includegraphics[width=0.45\textwidth]{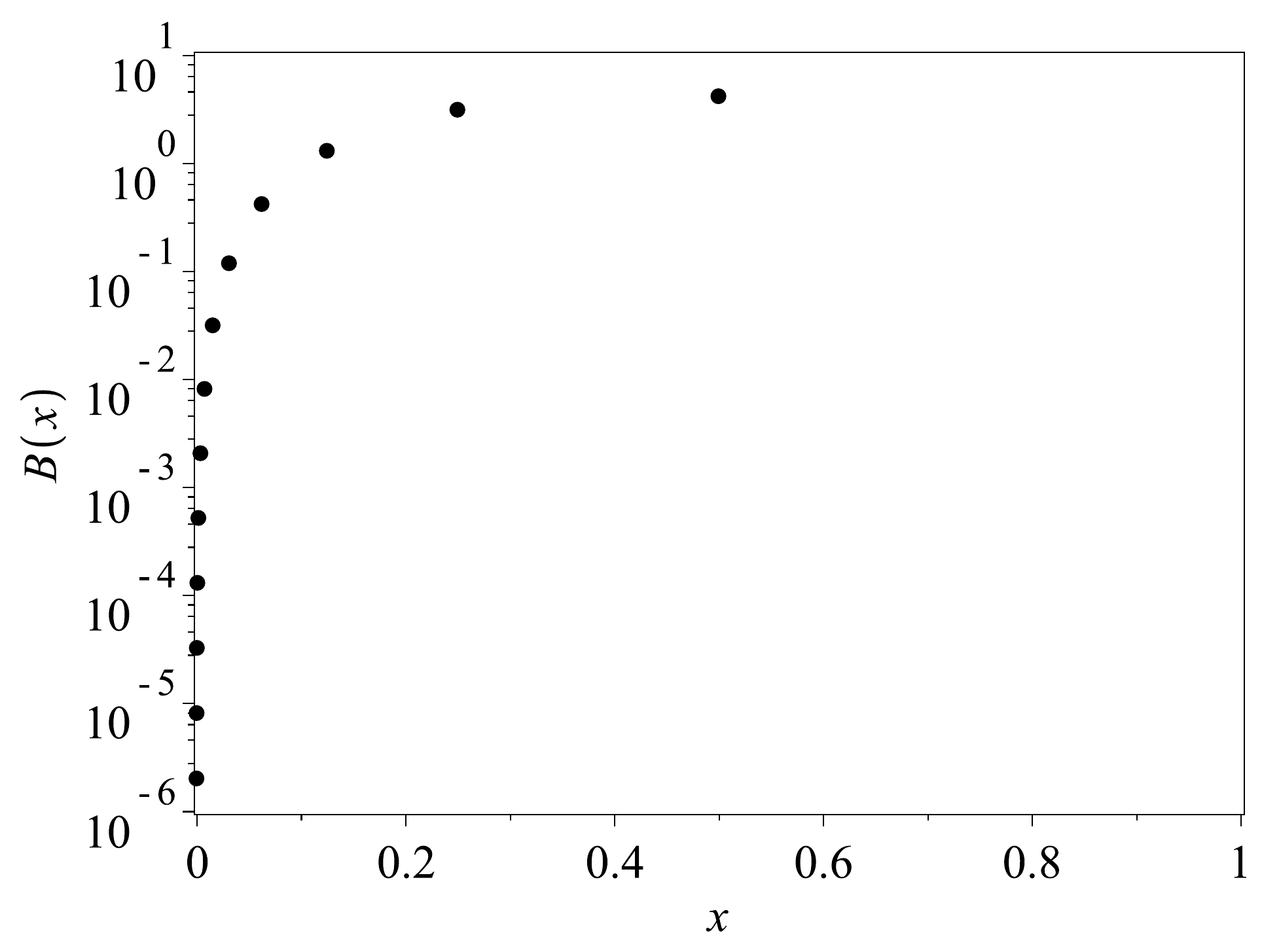}
    \caption{The condition number for Wilkionson's second example polynomial ($C_{20}$). In contrast to his first test problem, it is well-conditioned.}
    \label{Figure4c20}
\end{figure}

\begin{figure}[H]
    \centering
    \includegraphics[width=0.45\textwidth]{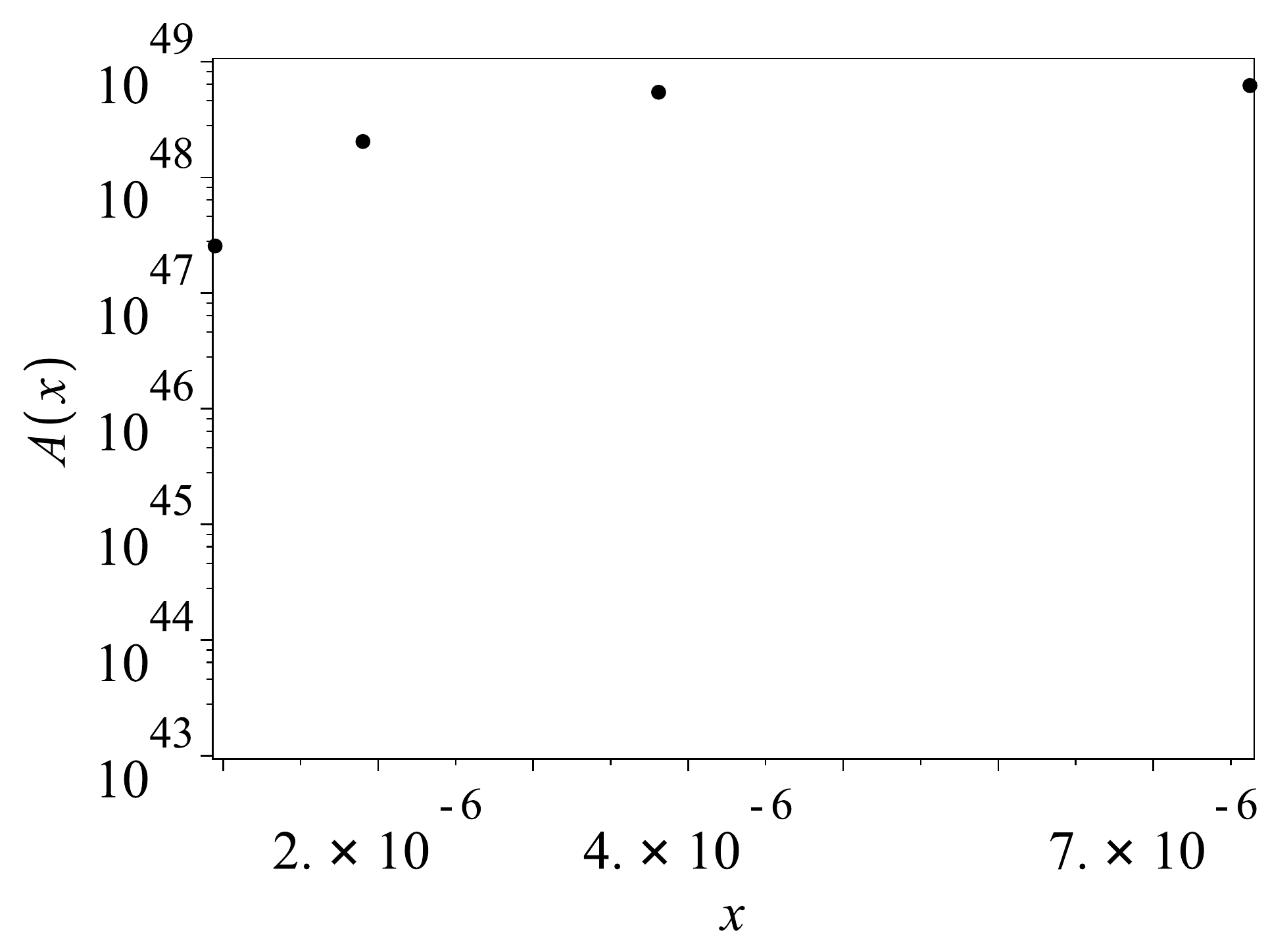}
    \caption{A portion of the condition number of $C_{20}$ in the Lagrange basis on the nodes $k/20$, $0\leq k\leq 20$.}
    \label{fig5}
\end{figure}

\begin{figure}[H]
    \centering
    \includegraphics[width=0.5\textwidth]{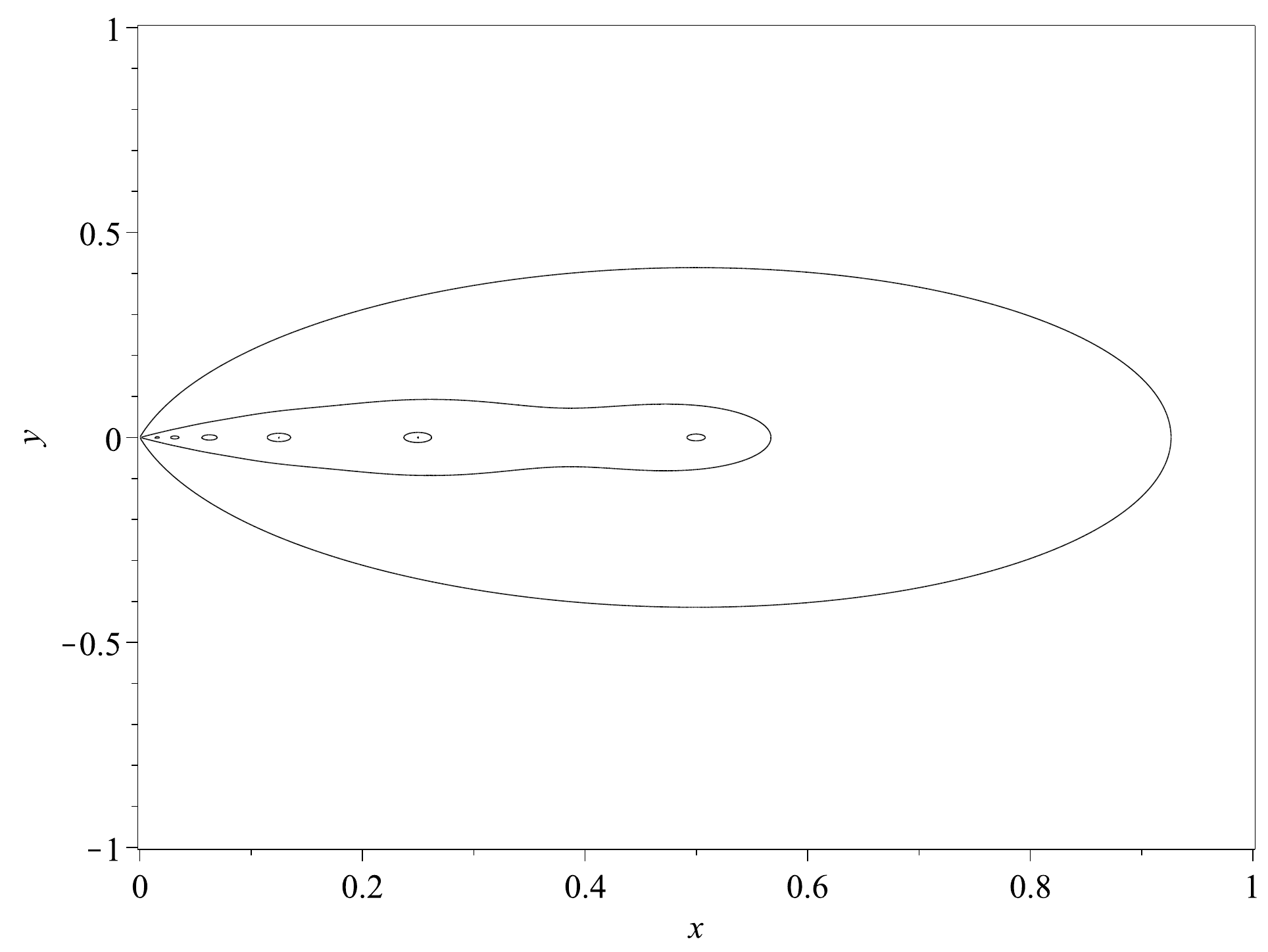}
    \caption{The pseudozeros of $C_{20}$. The contour levels are $10^{-1}$, $10^{-2}$, $10^{-3}$, $10^{-4}$, $10^{-6}$ and $10^{-8}$.}
    \label{pseudoC20}
\end{figure}

\begin{figure}[H]
    \centering
    \includegraphics[width=0.5\textwidth]{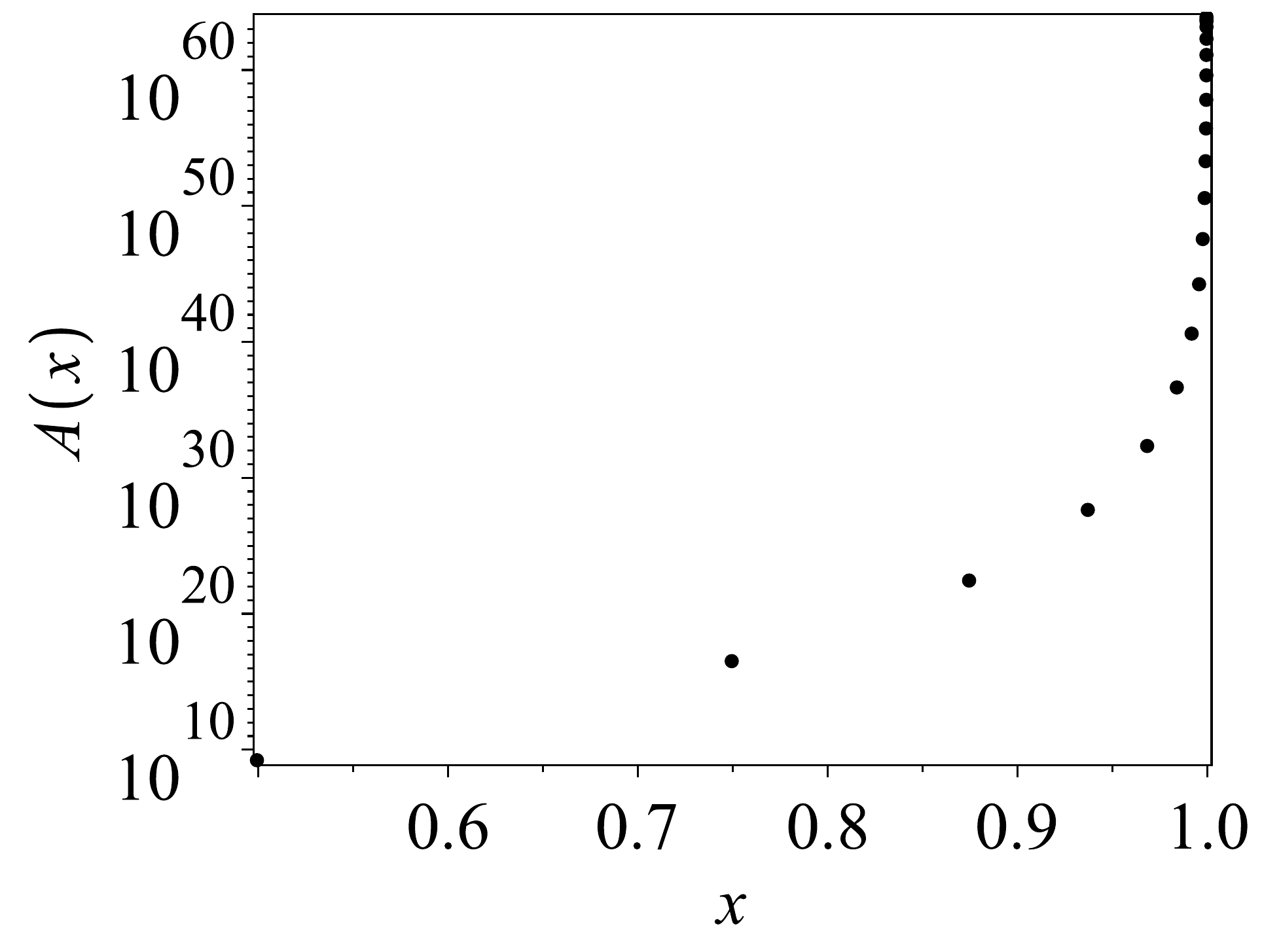}
    \caption{The condition number of $S_{20}$.}
    \label{FigureS20}
\end{figure}

\begin{figure}[H]
    \centering
    \includegraphics[width=0.5\textwidth]{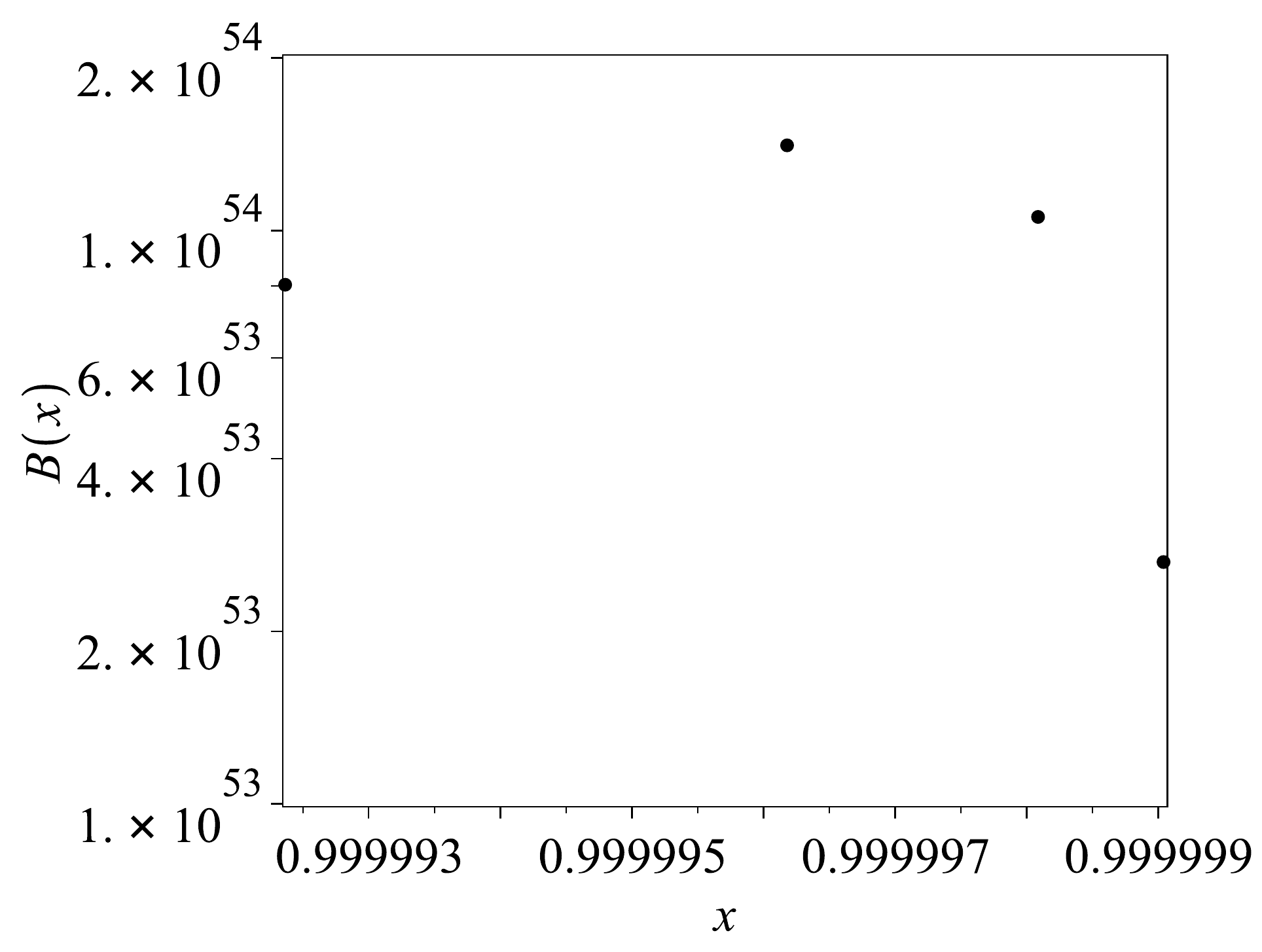}
    \caption{The condition number of $S_{20}$ in a Lagrange basis.}
    \label{figs20}
\end{figure}

\begin{figure}[H]
    \centering
    \includegraphics[width=0.5\textwidth]{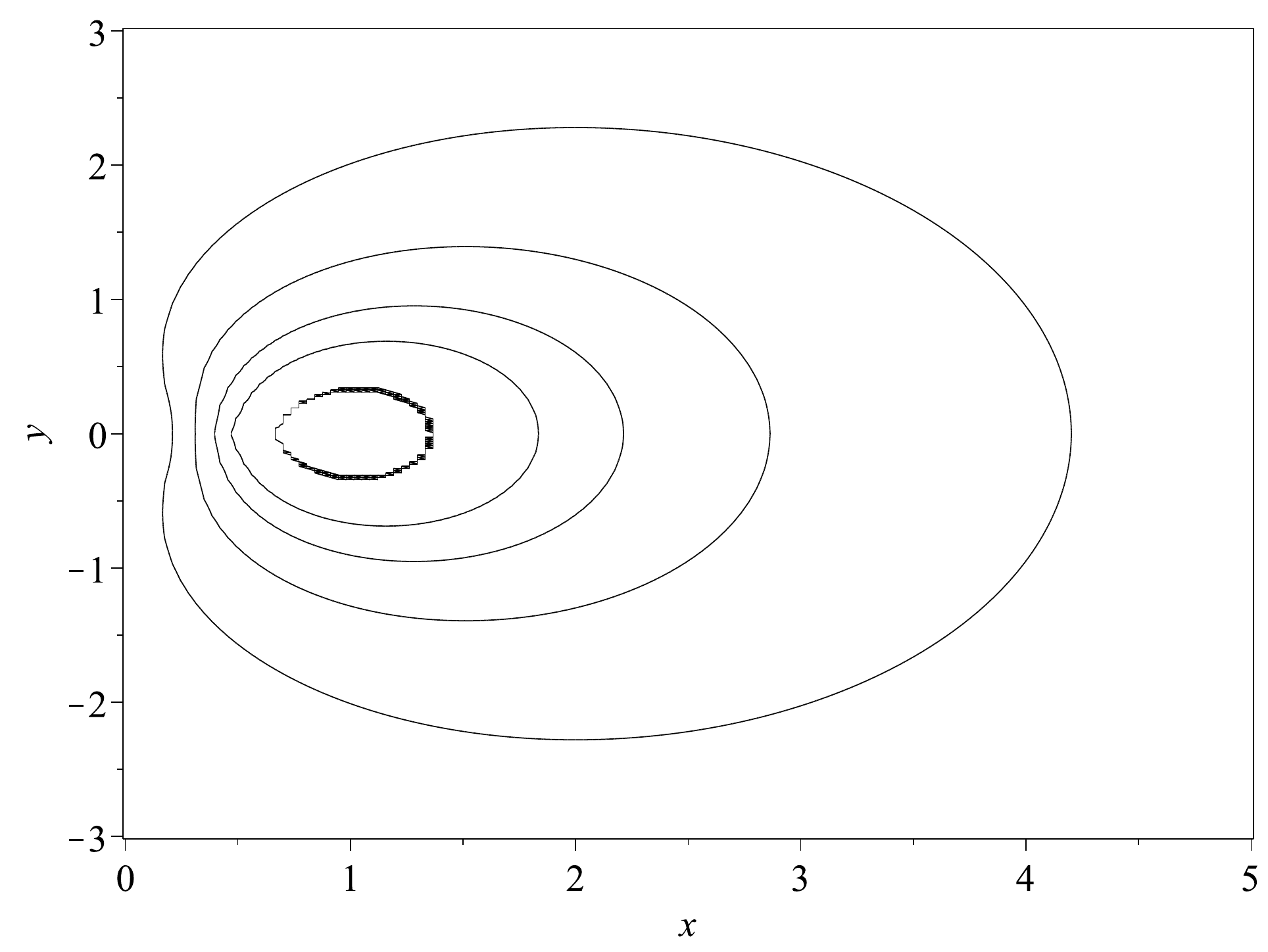}
    \caption{The pseudozeros of $S_{20}$. The contour levels are $10^{-4}$, $10^{-6}$, $10^{-8}$, $10^{-10}$ and $10^{-15}$.}
    \label{pseudoS20}
\end{figure}

The story of Wilkinson's second example is somehow more strange.
The polynomial is
\begin{equation}
    C_{20}(x) = \prod_{k=1}^{20} (x - 2^{-k})
\end{equation}
and the roots are $\sfrac{1}{2}, \sfrac{1}{4}, \sfrac{1}{8}, \sfrac{1}{16}, \ldots, \sfrac{1}{2^{20}}$.
Wilkinson expected that the clustering of roots near zero would cause difficulty for his rootfinder, once the polynomial was expanded:
\begin{equation}
    C_{20}(x) = x^{20} - \left ( \sum_{k=1}^{20} \frac{1}{2^k} \right )x^{19} + \cdots + \prod_{k=1}^{20}2^{-k} \>.
\end{equation}
But his program had no difficulty at all!
This is because the monomial basis is, in fact, quite well-conditioned near $x=0$, and the condition number for this polynomial can be seen in Figure \ref{Figure4c20} on $0 \le x \le 1$.
In contrast, the condition number for evaluation using the Lagrange basis on equally-spaced nodes in $[0,1]$, plus either $x_0 = 0$ or $x_0 = 1$, is horrible: for $N = 20$ it is already $10^{48}$, see Figure \ref{fig5}.
This computation conforms to Wilkinson's intuition that things can go wrong if roots are clustered.
Also we can see the pseudozeros of $C_{20}$ in Figure \ref{pseudoC20}. The required perturbations needed to make visible changes are quite large: these roots are not very sensitive to changes in the monomial basis coefficients.
Another way to see this is to look at a problem where the roots are clustered at $1$, not at $0$:
\begin{equation}
    S_{20} = \prod_{k=1}^{N} \left (x - (1 - 2^{-k}) \right )=\sum_{k=0}^N s_k x^k
\end{equation}
In this case the condition number is presented in Figure \ref{FigureS20}, and is huge. This polynomial is very sensitive to changes in the monomial basis coefficients.
The condition number for evaluation using a Lagrange basis for $S_{20}$ is shown in Figure \ref{figs20} (zoomed in for emphasis). Here the Lagrange basis is also very sensitive. 
Consider also the pseudozeros of $S_{20}$ in Figure \ref{pseudoS20}. The contour levels in Figure \ref{pseudoS20} are (from the outside in) $10^{-4}$, $10^{-6}$, $10^{-8}$, $10^{-10}$ and $10^{-15}$. In order to see a better view of the pseudozeros, let's consider the the first contour, $10^{-4}$, which is the biggest curve in Figure \ref{pseudoS20}. We know that 
\begin{equation}
S_{20}+\Delta S= \sum_{k=0}^{20}s_k(1+\delta_k)x^k
\end{equation}
where $\Delta S=s_0\delta_0+s_1\delta_1 x+s_2\delta_2 x^2+ \cdots + s_{20}\delta_{20} x^{20}$. Now if we choose a point between contour levels $10^{-4}$ and $10^{-6}$, for example $p=3-1.5i$, we can see that $p$ is a zero of some $S_{20}+\Delta S(x)$ with all coefficients of $\Delta S$ that have $|\delta_k| < 10^{-4}$. These are all small relative perturbations, that means everything inside the contour level $10^{-4}$ is a zero of a polynomial that is reasonably close to $S_{20}$. This is somehow backward error. So everything inside the contour level $10^{-4}$ is a zero of a polynomial closer to $S_{20}$ (in this sense) than $10^{-4}$. Everything inside the contour level $10^{-6}$ is a zero of a polynomial closer that $10^{-6}$ to $S_{20}$, and so on.
Notice that the innermost contour, corresponding to $10^{-15}$, is visible to the eye. This means that trivial (unit roundoff level in double precision) changes in the coefficients make visible changes in the root. 

\subsubsection*{Concluding remarks on the Wilkinson rootfinding examples}
The first example polynomial, $\prod_{k=1}^{20}(x-k)$, is nearly universally known as a surprising example. Yet there are very few places where one sees an elementary exposition of Wilkinson's theory of conditioning using this example, which is itself surprising because the theory was essentially born from this example. We have here illustrated Wilkinson's theory, as refined by Farouki and Rajan, for the students.\\
\begin{displayquote}
``For accidental historical reasons therefore backward error analysis is always introduced in connexion with matrix problems. In my opinion the ideas involved are much more readily absorbed if they are presented in connexion with polynomial equations. Perhaps the fairest comment would be that polynomial equations narrowly missed serving once again in their historical didactic role and rounding error analysis would have developed in a more satisfactory way if they had not.''
\\[5pt]
\rightline{{\rm --- James H. Wilkinson,~\cite{wilkinson}}}
\end{displayquote}
\subsubsection*{A final word for the instructor}
Backward error analysis is difficult at first for some kinds of students. The conceptual problem is that people are conditioned to think of mathematical problems as being exact; some indeed are, but many come from physical situations and are only models, with uncertain data. The success of BEA for floating point is to put rounding errors on the same footing as data or modeling errors, which have to be studied anyway. This is true even if the equations are solved exactly, by using computer algebra! The conditioning theory for polynomials discussed here allow this to be done quite flexibly, and are useful part of the analyst's repertoire. Students need to know this.

\bibliographystyle{plain}
\bibliography{bibliography}

\end{document}